\documentclass[12pt,leqno]{amsart}\usepackage[latin9]{inputenc}\usepackage{verbatim,amsthm,amstext,amssymb,color,amscd,bm}\makeatletter
\usepackage{hyperref}\hypersetup{colorlinks=true,citecolor=black,linkcolor=black}\usepackage[shortlabels]{enumitem}\numberwithin{equation}{section}
\theoremstyle{plain}\newtheorem{theorem}{Theorem}[section]\newtheorem{lemma}[theorem]{Lemma}
\theoremstyle{definition}\newtheorem{remark}[theorem]{Remark}\DeclareMathOperator{\Widim}{\mathrm{Widim}}\DeclareMathOperator{\mdim}{\mathrm{mdim}}
\makeatother
\begin{document}
\title[A minimal subshift of full mean dimension]{The Hilbert cube contains a minimal subshift\\of full mean dimension}
\author[Lei Jin]{Lei Jin}
\address{Lei Jin: School of Mathematics, Sun Yat-sen University, Guangzhou, China}
\email{jinleim@mail.ustc.edu.cn}
\author[Yixiao Qiao]{Yixiao Qiao}
\address{Yixiao Qiao (Corresponding author): School of Mathematics and Statistics, Guangdong University of Technology, Guangzhou, China}
\email{yxqiao@mail.ustc.edu.cn}
\begin{abstract}
We construct a minimal dynamical system of mean dimension equal to $1$, which can be embedded in the shift action on the Hilbert cube $[0,1]^\mathbb{Z}$. This clarifies a seemingly plausible impression about embedding possibility in relation to mean dimension. Our result finally leads to a full understanding of a pair of exact ranges of all the possible values of mean dimension, within which there will always be a minimal dynamical system that can/cannot be embedded in the shift action on the Hilbert cube.
\end{abstract}
\maketitle

\section{Main result}
Mean dimension is a topological invariant of dynamical systems, which originates with Misha Gromov \cite{Gromov} around 1999. It is closely connected with, and has been deeply applied to the \textit{embedding problem}, in particular, deciding if a minimal dynamical system can be embedded in the shift action on the Hilbert cube $[0,1]^\mathbb{Z}$. As follows is a brief review of the latest progress in this direction. With the necessary historical background we shall state our main theorem very quickly. All the precise notions and notations related to the statements can be found in Section 2, while a constructive proof of the main theorem is located in Section 3.

We denote by $([0,1]^\mathbb{Z},\sigma)$ the shift action on the Hilbert cube $[0,1]^\mathbb{Z}$ (note that its mean dimension is equal to $1$). Let $(X,T)$ be a \textit{minimal} dynamical system, whose mean dimension is denoted by $\mdim(X,T)$ which takes values within $[0,+\infty]$. The following statements are classically known:
\begin{itemize}
\item If $1<\mdim(X,T)\le+\infty$ then $(X,T)$ cannot be embedded in $([0,1]^\mathbb{Z},\sigma)$.
\item If $0\le\mdim(X,T)<1/2$ then $(X,T)$ can be embedded in $([0,1]^\mathbb{Z},\sigma)$.
\item If $1/2\le\mdim(X,T)<1$, then it is possible for $(X,T)$ to be embedded in $([0,1]^\mathbb{Z},\sigma)$, whereas it is also possible that $(X,T)$ cannot be embedded in $([0,1]^\mathbb{Z},\sigma)$.
\item Suppose $\mdim(X,T)=1$. It may happen that $(X,T)$ cannot be embedded in $([0,1]^\mathbb{Z},\sigma)$.
\end{itemize}
Here let us make some remarks on the above assertions. The first implication follows directly from definition. The second significant result is due to Gutman and Tsukamoto \cite{GT} (for a corresponding result of $\mathbb{Z}^k$-actions we refer to \cite{GQT}). In relation to the third and fourth cases, embeddable and non-embeddable examples were constructed by Lindenstrauss--Weiss \cite{LW} and Lindenstrauss--Tsukamoto \cite{LT}, respectively. To be precise, Lindenstrauss and Tsukamoto \cite{LT} showed that for any given $r\ge1/2$ there is a minimal dynamical system of mean dimension equal to $r$, which cannot be embedded in $([0,1]^\mathbb{Z},\sigma)$; Lindenstrauss and Weiss \cite{LW} proved that for any given $0\le r<1$ there exists a minimal dynamical system of mean dimension equal to $r$, which can be embedded in $([0,1]^\mathbb{Z},\sigma)$.

A careful reader may observe that in order to have a complete picture of the situation of embedding minimal dynamical systems in the Hilbert cube, there is only one issue that remains open, namely, the embedding possibility at the critical value $1$ (for mean dimension) is not clear to us. Formally, we study the problem as follows:
\begin{itemize}
\item Assume that a minimal dynamical system can be embedded in $([0,1]^\mathbb{Z},\sigma)$. Does this imply that it has mean dimension strictly less than $1$?
\end{itemize}
In other words, this problem asks if there exists a minimal dynamical system of mean dimension equal to $1$, which can be embedded in $([0,1]^\mathbb{Z},\sigma)$. The aim of this paper is to solve this problem.

In fact, we observe the following three aspects:
\begin{enumerate}
\item With a view towards mean dimension theory: The approach of any previously mentioned results \cite{LW,LT} reveals that if we try to construct a minimal dynamical system of mean dimension equal to $1$, then a construction within the framework of $[0,1]^\mathbb{Z}$ is \textit{not} adequate for our purpose, and instead, we need consider, with the same method, the shift action on $([0,1]^2)^\mathbb{Z}$, which has mean dimension equal to $2$, rather than consider the shift action on $[0,1]^\mathbb{Z}$.
\item With a view towards entropy theory: We note moreover that a highly similar circumstance took place in topological entropy in connection with the shift action over any alphabet of finite cardinality; strictly speaking, such an ``analogue'' asserts that any \textit{proper} subsystem of the shift action on $\{k\in\mathbb{Z}:1\le k\le N\}^\mathbb{Z}$ (where $N$ is a positive integer) does \textit{not} have full topological entropy (i.e. its topological entropy must be strictly less than $\log N$).
\item With a view towards dimension theory: A celebrated theorem in dimension theory states that if a compact metrizable space of (topological) dimension $D$ (where $D$ is a positive integer) can be (topologically) embedded in $[0,1]^D$, then it must contain (a copy of) $[0,1]^D$ (i.e. $[0,1]^D$ can be embedded in it).
\end{enumerate}

Apparently, all of these three notable phenomena lead naturally to a seemingly plausible impression, i.e. it would be true that minimal dynamical systems satisfying the embeddable property are \textit{not} able to possess \textit{full} mean dimension. However, this impression turns out to be \textit{false}. Somewhat surprisingly, we do construct a minimal dynamical system successfully, which satisfies both embeddable and full mean dimensional conditions. This is our main result.

\begin{theorem}[Main theorem]\label{main}
There exists a minimal dynamical system of mean dimension equal to $1$, which can be embedded in the shift action on the Hilbert cube $[0,1]^\mathbb{Z}$.
\end{theorem}

In contrast to the seemingly reasonable observation as stated above, Theorem \ref{main} enables a different (and unexpected) behaviour to become clarified. Furthermore, our main theorem, together with the previous remarkable results, eventually allows the embeddability of minimal dynamical systems in the Hilbert cube to be fully understood with a view towards mean dimension. More precisely, we have now gotten a full understanding of (a pair of) the \textit{exact} ranges (I) (i.e. $[0,1]$) and (II) (i.e. $[1/2,+\infty]$) of \textit{all the possible values} of mean dimension, such that for any value assigned within (I) (resp. within (II)) there always exists a minimal dynamical system of mean dimension equal to this appointed value, which can be (resp. cannot be) embedded in the shift action on the Hilbert cube. Notice that with Theorem \ref{main} the range (I) finally becomes \textit{exact}, whereas the range (II) is \textit{exact} because of \cite{LT} for the \textit{existence} of such a system of mean dimension within $[1/2,+\infty]$ and because of \cite{GT} for the \textit{non-existence} of such a system of mean dimension within $[0,1/2)$.

The strategy we shall adopt generally follows the technical framework of the ``block-type'' induction. Nevertheless, we have to make each step in our construction sufficiently delicate. The key ingredient of our idea is to produce a dense subset of the alphabet $[0,1]$ ``more gently'' (i.e. with a sequence of closed intervals instead of points).

\begin{remark}
The same statement (as in Theorem \ref{main}) also applies to the alphabet $[0,1]^D$; namely, the following assertion is true.
\begin{itemize}\item
For any positive integer $D$ (possibly $+\infty$) there is a minimal dynamical system of mean dimension $D$, which can be embedded in the shift action on $([0,1]^D)^\mathbb{Z}$.
\end{itemize}
Meanwhile, the alphabet $[0,1]^D$ of topological dimension $D$ may be replaced by some other compact metrizable spaces provided they have some \textit{nice} structure, e.g. polyhedrons $P$ of topological dimension $\dim(P)$. Besides, we notice that it is straightforward to generalize our main theorem to actions of countably infinite amenable groups (with almost no additional effort while using material on tilings of amenable groups).
\end{remark}

\subsection*{Acknowledgements}
The initial proof we provided in a previous version of this paper contained a gap. Professor Masaki Tsukamoto pointed out this to us, and meanwhile, he also kindly explained to us how to fix it. We would like to thank him for his warm help. We are also grateful to the anonymous referee for his/her insightful comments and helpful suggestions which improve this paper greatly. L. Jin was supported by NNSF of China No. 12201653. Y. Qiao was supported by NNSF of China No. 12371190.

\section{Terminologies}
This section is devoted to a brief review of all the fundamental definitions appearing in the previous section. By a (topological) \textbf{dynamical system} we shall understand a pair $(X,T)$, where $X$ is a compact metrizable space and $T:X\to X$ is a homeomorphism. An important class of dynamical systems are minimal systems. A dynamical system $(X,T)$ is said to be \textbf{minimal} if for every point $x\in X$ the set $\{T^nx:n\in\mathbb{Z}\}$ is dense in $X$. Among other typical examples of dynamical systems, probably the most canonical ones are the \textbf{shift action} on the Hilbert cubes $([0,1]^D)^\mathbb{Z}$ (where $D$ is a positive integer or $+\infty$), which we denote by $(([0,1]^D)^\mathbb{Z},\sigma)$, defined as follows: $$\sigma:([0,1]^D)^\mathbb{Z}\to([0,1]^D)^\mathbb{Z},\quad(x_n)_{n\in\mathbb{Z}}\mapsto(x_{n+1})_{n\in\mathbb{Z}}.$$ For simplicity we assume $D=1$. We say that a dynamical system $(X,T)$ can be \textbf{embedded} in $([0,1]^\mathbb{Z},\sigma)$ if there is an equivariant topological embedding $f:X\to[0,1]^\mathbb{Z}$, namely, a homeomorphism of $X$ into $[0,1]^\mathbb{Z}$ satisfying $f\circ T=\sigma\circ f$. Such a mapping $f$ is usually called an embedding of the dynamical system $(X,T)$ in $([0,1]^\mathbb{Z},\sigma)$. As we have seen in Section 1, this paper mainly focuses on the embeddability of minimal dynamical systems in $([0,1]^\mathbb{Z},\sigma)$.

We denote by $\dim(P)$ the topological dimension (i.e. the Lebesgue covering dimension) of a compact metrizable space $P$ (which is always assumed to be nonempty). Let $X$ and $P$ be two compact metrizable spaces and $\rho$ a compatible metric on $X$. For $\epsilon>0$ a continuous mapping $f:X\to P$ is called an \textbf{$\epsilon$-embedding} with respect to $\rho$ if $f(x)=f(x^\prime)$ implies $\rho(x,x^\prime)<\epsilon$, for all $x,x^\prime\in X$. Let $\Widim_\epsilon(X,\rho)$ be the minimum topological dimension $\dim(P)$ of a compact metrizable space $P$ which admits an $\epsilon$-embedding $f:X\to P$ with respect to $\rho$.
\begin{remark}
We may verify that the topological dimension of $X$ may be recovered by $\dim(X)=\lim_{\epsilon\to0}\Widim_\epsilon(X,\rho)$.
\end{remark}
Let $(X,T)$ be a dynamical system with a compatible metric $\rho$ on $X$. For every positive integer $n$ we define on $X$ a compatible metric $\rho_n$ as follows: $$\rho_n(x,x^\prime)=\max_{0\le i<n}\rho(T^ix,T^ix^\prime),\quad\forall\,x,x^\prime\in X.$$ The \textbf{mean dimension} of $(X,T)$ is defined by $$\mdim(X,T)=\lim_{\epsilon\to0}\lim_{n\to+\infty}\frac{\Widim_\epsilon(X,\rho_n)}{n}.$$ It is well known that the limits in the above definition always exist, and the value $\mdim(X,T)$ is independent of the choices of a compatible metric $\rho$ on $X$ and a F\o lner sequence $\{F_n\}_{n=1}^{+\infty}$ (instead of $\{k\in\mathbb{Z}:0\le k<n\}_{n=1}^{+\infty}$) of $\mathbb{Z}$.

We note that the mean dimension of $(([0,1]^D)^\mathbb{Z},\sigma)$ is equal to $D$, where $D$ is a positive integer or $+\infty$. For details we refer the reader to \cite{LW}. Clearly, if a dynamical system can be embedded in $(([0,1]^D)^\mathbb{Z},\sigma)$, then its mean dimension must be less than or equal to $D$.

\section{Proof of Theorem \ref{main}}
\subsection{Construction of $(X,\sigma)$}
We are going to construct a minimal dynamical system $(X,\sigma)$ such that $X$ is a closed and shift-invariant subset of the Hilbert cube $[0,1]^\mathbb{Z}$ and that the mean dimension of $(X,\sigma)$ is equal to $1$.

First of all, we fix a compatible metric $d$ on $[0,1]^\mathbb{Z}$ as follows: $$d(x,x^\prime)=\sum_{n\in\mathbb{Z}}\frac{|x_n-x^\prime_n|}{2^{|n|}},\quad\left(x=(x_n)_{n\in\mathbb{Z}},x^\prime=(x^\prime_n)_{n\in\mathbb{Z}}\in[0,1]^\mathbb{Z}\right).$$ For $x=(x_n)_{n\in\mathbb{Z}}\in[0,1]^\mathbb{Z}$ and two integers $t\le t^\prime$ we set $x|_t^{t^\prime}=(x_n)_{t\le n\le t^\prime}$. For a closed interval $I$ we denote its length by $|I|$.

The construction will be fulfilled by induction. To start with, let us explain the intuitive meaning of our notations shortly. For a nonnegative integer $k$ we build in the $k$-th step a closed shift-invariant subset $X_k$ of $[0,1]^\mathbb{Z}$, which is generated from a ``block'' $B_k$ of the form $\prod_{j=1}^{b_k}I_j^{(k)}$, where each $I_j^{(k)}$ is a closed subinterval of $[0,1]$ and where the positive integer $b_k$ indicates the ``length'' of the block $B_k$. Roughly speaking, employing a family of closed intervals in our construction is actually to increase dimension. But this is \textit{not} enough for its mean dimension to be dominated from below. Therefore we are going through the following approach: When we deal with the next block $B_{k+1}$ in the $(k+1)$-th step we have to ``copy'' the block $B_k$ sufficiently many times in order to occupy a rather large proportion of positions of the block $B_{k+1}$. The positive integer $r_k$ is to describe the number of copies of the block $B_k$ in the block $B_{k+1}$, whereas for an integer $0\le s<k+1$, the number $\eta(s,k+1)>0$ is to dominate the ``density'' of those intervals of length greater than or equal to $1/2^{s(s+1)/2}$ in the block $B_{k+1}$, which will strictly decrease as $k$ increases, while the number $\eta(s)>0$ is to guarantee that the lower bound for the value $\eta(s,k+1)$ is always under control as $k+1$ ranges over all those integers greater than $s$, which will finally converge to $1$ as $s$ goes to $+\infty$.

\textbf{Step} $-\infty$. We take a two-parameter sequence $\{\eta(s,k):0\le s<k;s,k\in\mathbb{Z}\}$ of real numbers and a one-parameter sequence $\{\eta(s)\}_{s=0}^{+\infty}$ of real numbers satisfying all the following conditions (the existence of such two sequences is obvious):
\begin{enumerate}
\item For any two integers $0\le s<k$: $\;0<\eta(s,k)<1$; $\,0<\eta(s)<1$.
\item For any nonnegative integer $s$ the one-parameter sequence $\{\eta(s,k)\}_{k=s+1}^{+\infty}$ is \textit{strictly} decreasing, namely $$\eta(s,k+1)<\eta(s,k),\;\quad\forall\,0\le s<k\quad(s,k\in\mathbb{Z}).$$
\item For any nonnegative integer $s$ the one-parameter sequence $\{\eta(s,k)\}_{k=s+1}^{+\infty}$ is bounded by $\eta(s)$ from below, i.e. $$\eta(s,k)>\eta(s),\;\quad\forall\,0\le s<k\quad(s,k\in\mathbb{Z}).$$
\item $\eta(s)\to1$ as $s\to+\infty$.
\end{enumerate}
We fix them throughout this section. We would like to remind the reader that we shall use these conditions \textit{implicitly} in the sequel.

\textbf{Step} $0$. We take $b_0=1$, $B_0=I_1^{(0)}=[0,1]$ and $X_0=[0,1]^\mathbb{Z}$.

\textbf{Step} $1$. We divide $I_1^{(0)}=[0,1]$ equally into $2$ closed subintervals of length $1/2$ as follows: $$I_{1,1}^{(1)}=[0,1/2],\quad\;I_{1,2}^{(1)}=[1/2,1].$$ We take a positive integer $r_0$ with $$\frac{r_0}{r_0+2}\ge\eta(0,1).$$ We put $b_1=r_0+2$ and set $$B_1=(B_0)^{r_0}\times I_{1,1}^{(1)}\times I_{1,2}^{(1)}\subset[0,1]^{b_1}.$$ We rename the $b_1$ closed subintervals of $[0,1]$ appearing in the above product with indices, and rewrite $B_1$ as follows: $B_1=\prod_{j=1}^{b_1}I_j^{(1)}$. More precisely, we notice here (in relation to the indices) that $B_0=I_1^{(1)}$ and that the sets $I_{b_1-1}^{(1)}$ and $I_{b_1}^{(1)}$ form a cover of $B_0$. Let $X_1$ be the set of all those $x=(x_n)_{n\in\mathbb{Z}}\in[0,1]^\mathbb{Z}$ satisfying the following condition:
\begin{itemize}\item
There exists some integer $0\le l\le b_1-1$ such that $x|_{l+b_1m+1}^{l+b_1(m+1)}\in B_1$ for all $m\in\mathbb{Z}$.
\end{itemize}
Obviously, $X_1$ is a nonempty closed and shift-invariant subset of $X_0=[0,1]^\mathbb{Z}$.

To proceed, we assume for a positive integer $k$ that $b_{k-1}$, $B_{k-1}$ and $X_{k-1}$ have already been generated in Step $(k-1)$. We now generate $b_k$, $B_k$ and $X_k$.

\textbf{Step} $k$. We write $B_{k-1}=\prod_{j=1}^{b_{k-1}}I_j^{(k-1)}$, where every $I_j^{(k-1)}$ is a closed subinterval of $[0,1]$. For each $1\le j\le b_{k-1}$ we divide $I_j^{(k-1)}$ equally into $2^k$ closed subintervals of length equal to $|I_j^{(k-1)}|/2^k$, which are denoted by $I_{j,i}^{(k)}$ with $1\le i\le 2^k$. We take a positive integer $r_{k-1}$ sufficiently large such that for every integer $0\le s\le k-1$ $$\frac{\left(\#\{1\le j\le b_{k-1}:|I_j^{(k-1)}|\ge1/2^{s(s+1)/2}\}\right)\cdot r_{k-1}}{b_{k-1}\cdot(r_{k-1}+2^{kb_{k-1}})}\ge\eta(s,k),$$ where the symbol $\#$ records the cardinality of a set. We put $b_k=b_{k-1}\cdot(r_{k-1}+2^{kb_{k-1}})$ and set $$B_k\,=\,(B_{k-1})^{r_{k-1}}\,\times\,\prod_{i_1,\dots,i_{b_{k-1}}\in\{1,\dots,2^k\}}\;\prod_{j=1}^{b_{k-1}}\,I_{j,i_j}^{(k)}\;\,\subset\,(B_{k-1})^{r_{k-1}+2^{kb_{k-1}}}\,.$$ We note that $B_k$ is a closed subset of $[0,1]^{b_k}$. We rewrite $B_k$ as $B_k=\prod_{j=1}^{b_k}I_j^{(k)}$, where $I_1^{(k)},\dots,I_{b_k}^{(k)}$ are closed subintervals of $[0,1]$. To be precise here we remark that in relation to the indices we may ensure in addition: $$\bigcup_{m=r_{k-1}}^{r_{k-1}+2^{kb_{k-1}}-1}\,\prod_{j=b_{k-1}m+1}^{b_{k-1}(m+1)}I_j^{(k)}\,=\,B_{k-1}\,=\,\prod_{j=1}^{b_{k-1}}I_j^{(k)}.$$ We let $X_k$ be the set consisting of all the points $x=(x_n)_{n\in\mathbb{Z}}\in[0,1]^\mathbb{Z}$ satisfying the following condition:
\begin{itemize}\item
There is some integer $0\le l\le b_k-1$ such that $x|_{l+b_km+1}^{l+b_k(m+1)}\in B_k$ for all $m\in\mathbb{Z}$.
\end{itemize}
It follows that $X_k$ is a (nonempty) closed and shift-invariant subset of $X_{k-1}\subset[0,1]^\mathbb{Z}$.

\textbf{Step} $+\infty$. The induction has now been completed, as we have already generated $b_k$, $B_k$ and $X_k$ for all nonnegative integers $k$. To end the construction, we finally take the intersection as follows: $$X=\bigcap_{k=0}^{+\infty}X_k.$$ Since $\{X_k\}_{k=0}^{+\infty}$ is a decreasing sequence of nonempty closed shift-invariant subsets of $[0,1]^\mathbb{Z}$, $X$ is a nonempty closed shift-invariant subset of $[0,1]^\mathbb{Z}$ as well. Thus, $(X,\sigma)$ becomes a dynamical system. In what follows we need verify that it does satisfy the required conditions.

\subsection{Minimality of $(X,\sigma)$}
We will show that the dynamical system $(X,\sigma)$ is minimal. In fact, it suffices to show that for any $x,y\in X$ and any $\epsilon>0$ there exists some $M\in\mathbb{Z}$ such that $d(\sigma^M(x),y)<\epsilon$. Let us fix $x=(x_n)_{n\in\mathbb{Z}},y=(y_n)_{n\in\mathbb{Z}}\in X$ and $\epsilon>0$.

For each positive integer $k$ we set $$M_k=\max\left\{|I_j^{(k)}|:b_{k-1}\cdot r_{k-1}+1\le j\le b_k\right\}.$$ This sequence $\{M_k\}_{k=1}^{+\infty}$ is deceasing and converges to $0$ as $k$ goes to $+\infty$. Actually, it follows from the construction that the value $M_k$ is equal to $1/2^k$.

We choose a positive integer $L$ depending only on $\epsilon>0$ such that if two points $x^\prime=(x^\prime_n)_{n\in\mathbb{Z}}$ and $x^{\prime\prime}=(x^{\prime\prime}_n)_{n\in\mathbb{Z}}$ coming from $[0,1]^\mathbb{Z}$ satisfy that $|x^\prime_n-x^{\prime\prime}_n|<\epsilon/2$ for any $-L\le n\le L$, then they will satisfy $d(x^\prime,x^{\prime\prime})<\epsilon$.

We take an integer $N>2L+1$ with $1/2^N<\epsilon/2$. Let us look at the $N$-th, $(N+1)$-th and $(N-1)$-th steps of the construction.

Since $x,y\in X\subset X_N$, there must exist two integers $0\le p,q\le b_N-1$ such that $$x|_{p+b_Nm+1}^{p+b_N(m+1)},y|_{q+b_Nm+1}^{q+b_N(m+1)}\in B_N,\quad\forall m\in\mathbb{Z}.$$ There are two integers $L_0\in[-L,L]$ and $m_0\in\{-1,0\}$, both of which are uniquely determined, satisfying that $$[-L,L_0-1]\subset[q+b_N(m_0-1)+1,q+b_Nm_0],$$$$[L_0,L]\subset[q+b_Nm_0+1,q+b_N(m_0+1)].$$ Here we assume by convention $[-L,-L-1]=\emptyset$.

Let us make it clearer with a short remark. Strictly speaking, the general case is where $L_0\in[-L+1,L]$ (and hence where we will have $L_0=q+b_Nm_0+1$), with an exception (which turns out to be simpler than the general case) if such an integer $L_0\in[-L+1,L]$ does not exist, for which we set $L_0=-L$. As we will see in a moment, the exceptional case is contained in (a part of) the general case. Therefore we assume without loss of generality $-L+1\le L_0\le L$.

Since $$y|_{q+b_Nm_0+1}^{q+b_N(m_0+1)}\in B_N=\bigcup_{m=r_N}^{r_N+2^{(N+1)b_N}-1}\prod_{j=b_Nm+1}^{b_N(m+1)}I_j^{(N+1)},$$ there is an integer $r_N\le m_1\le r_N+2^{(N+1)b_N}-1$ with $$y|_{q+b_Nm_0+1}^{q+b_N(m_0+1)}\in\prod_{j=b_Nm_1+1}^{b_N(m_1+1)}I_j^{(N+1)}.$$ Since $x\in X\subset X_{N+1}$ and since $$B_{N+1}=(B_N)^{r_N}\times\prod_{j=b_Nr_N+1}^{b_{N+1}}I_j^{(N+1)},$$ there is some integer $m_2$ depending on $m_1$ such that $$x|_{p+b_Nm_2+1}^{p+b_N(m_2+1)}\in\prod_{j=b_Nm_1+1}^{b_N(m_1+1)}I_j^{(N+1)}.$$ Thus, for every integer $1\le i\le b_N$ $$|x_{p+b_Nm_2+i}-y_{q+b_Nm_0+i}|\le M_{N+1}.$$ This implies that for every integer $n\in[q+b_Nm_0+1,q+b_N(m_0+1)]$, in particular, for every integer $n\in[L_0,L]$ $$|\sigma^{b_N(m_2-m_0)+p-q}(x)_n-y_n|\le M_{N+1}=1/2^{N+1}<\epsilon/2.$$

Since $$\sigma^{b_N(m_2-m_0)+p-q}(x)|_{q+b_N(m_0-1)+1}^{q+b_Nm_0},y|_{q+b_N(m_0-1)+1}^{q+b_Nm_0}\in B_N$$ and since $$2L+1\le N-1\le b_{N-1},$$ we deduce that for every integer $n\in[-L,L_0-1]$ $$|\sigma^{b_N(m_2-m_0)+p-q}(x)_n-y_n|\le M_N=1/2^N<\epsilon/2.$$

Thus, we conclude that for all integers $n\in[-L,L]$ $$|\sigma^{b_N(m_2-m_0)+p-q}(x)_n-y_n|<\epsilon/2,$$ which implies $$d(\sigma^{b_N(m_2-m_0)+p-q}(x),y)<\epsilon.$$

\subsection{Mean dimension of $(X,\sigma)$}
We shall show that $\mdim(X,\sigma)=1$. Since it is obvious that the dynamical system $(X,\sigma)$ has been embedded in $([0,1]^\mathbb{Z},\sigma)$ already, its mean dimension does not exceed $1$. Thus, it suffices to prove $\mdim(X,\sigma)\ge1$. This will end the paper.

We borrow a practical lemma \cite[Lemma 1.1.1]{Gromov}\cite[Example 2.1]{LT} as follows. The point which will be crucial to our argument is that the equality presented in this lemma is true for all positive integers $n$ and all sufficiently small $\delta>0$ (which should be \textit{strictly} less than $\tau>0$).
\begin{lemma}\label{widimcube}
For any $0<\delta<\tau$ and any positive integer $n$ $$\Widim_\delta\left([0,\tau]^n,d_{l^\infty}\right)=n.$$ Here $d_{l^\infty}$ is the compatible metric on $[0,\tau]^n$ defined by $$d_{l^\infty}\left((x_i)_{i=1}^n,(x^\prime_i)_{i=1}^n\right)=\max_{1\le i\le n}|x_i-x^\prime_i|,\;\quad\forall\,(x_i)_{i=1}^n,(x^\prime_i)_{i=1}^n\in[0,\tau]^n.$$
\end{lemma}

We fix a point $z\in X$ satisfying that $z|_{b_km}^{b_k(m+1)-1}\in B_k$ for any nonnegative integer $k$ and any $m\in\mathbb{Z}$. For each nonnegative integer $k$ we define a mapping as follows: $$F_k:B_k\to X,\,\quad x=(x_n)_{n=0}^{b_k-1}\,\mapsto\,F_k(x)=(F_k(x)_n)_{n\in\mathbb{Z}},$$$$F_k(x)_n=\begin{cases}x_n,\quad&0\le n\le b_k-1\\z_n,\quad&n\in\mathbb{Z}\setminus[0,b_k-1]\end{cases}.$$ It is clear that $F_k(x)$ (where $x\in B_k$) is indeed in $X$. Moreover, the mapping $F_k:B_k\to X$ is continuous and distance-increasing with respect to the compatible metric $d_{l^\infty}$ on $B_k$ and the compatible metric $d_{b_k}$ on $X$, i.e. it satisfies that $$d_{l^\infty}(x,x^\prime)\le d_{b_k}(F_k(x),F_k(x^\prime)),\quad\forall\,x,x^\prime\in B_k.$$ It follows that for every $\epsilon>0$ and every nonnegative integer $k$ $$\Widim_\epsilon(B_k,d_{l^\infty})\le\Widim_\epsilon(X,d_{b_k}).$$ Now let us take $\epsilon>0$ arbitrarily. We choose a nonnegative integer $s=s(\epsilon)$ depending only on $\epsilon>0$ such that $$1/2^{(s+1)(s+2)/2}\le\epsilon<1/2^{s(s+1)/2}.$$ Note that for all integers $k>s$ $$\frac{\#\{1\le j\le b_k:|I_j^{(k)}|\ge1/2^{s(s+1)/2}\}}{b_k}\ge\eta(s,k).$$ By Lemma \ref{widimcube} this implies that for all integers $k>s$ $$\frac{\Widim_\epsilon(B_k,d_{l^\infty})}{b_k}\ge\eta(s,k)\ge\eta(s).$$ Thus, for all integers $k>s$ $$\frac{\Widim_\epsilon(X,d_{b_k})}{b_k}\ge\eta(s).$$ Since $\epsilon>0$ is arbitrary and since $s\to+\infty$ as $\epsilon\to0$, we conclude with $$\mdim(X,\sigma)=\lim_{\epsilon\to0}\lim_{k\to+\infty}\frac{\Widim_\epsilon(X,d_{b_k})}{b_k}\ge\lim_{\epsilon\to0}\eta(s)=1.$$

\medskip

\end{document}